\documentclass[12pt]{article}

\usepackage{amsmath}
\usepackage{amssymb}

\def\hybrid{\topmargin 0pt      \oddsidemargin 0pt
        \headheight 0pt \headsep 0pt
        \voffset=-0.5cm
        \hoffset=-0.25in
        \textwidth 6.75in
        \textheight 9.5in       
        \marginparwidth 0.0in
        \parskip 5pt plus 1pt   \jot = 1.5ex}
\catcode`\@=11
\def\marginnote#1{}

\newcount\hour
\newcount\minute
\newtoks\amorpm
\hour=\time\divide\hour by60 \minute=\time{\multiply\hour by60
\global\advance\minute by-\hour}
\edef\standardtime{{\ifnum\hour<12 \global\amorpm={am}%
        \else\global\amorpm={pm}\advance\hour by-12 \fi
        \ifnum\hour=0 \hour=12 \fi
        \number\hour:\ifnum\minute<10 0\fi\number\minute\the\amorpm}}
\edef\militarytime{\number\hour:\ifnum\minute<10 0\fi\number\minute}

\def\draftlabel#1{{\@bsphack\if@filesw {\let\thepage\relax
   \xdef\@gtempa{\write\@auxout{\string
      \newlabel{#1}{{\@currentlabel}{\thepage}}}}}\@gtempa
   \if@nobreak \ifvmode\nobreak\fi\fi\fi\@esphack}
        \gdef\@eqnlabel{#1}}
\def\@eqnlabel{}
\def\@vacuum{}
\def\draftmarginnote#1{\marginpar{\raggedright\scriptsize\tt#1}}
\def\draftlabel#1{{\@bsphack\if@filesw {\let\thepage\relax
   \xdef\@gtempa{\write\@auxout{\string
      \newlabel{#1}{{\@currentlabel}{\thepage}}}}}\@gtempa
   \if@nobreak \ifvmode\nobreak\fi\fi\fi\@esphack}
        \gdef\@eqnlabel{#1}}
\def\@eqnlabel{}
\def\@vacuum{}
\def\draftmarginnote#1{\marginpar{\raggedright\scriptsize\tt#1}}

\def\draft{\oddsidemargin -.5truein
        \def\@oddfoot{\sl preliminary draft \hfil
        \rm\thepage\hfil\sl\today\quad\militarytime}
        \let\@evenfoot\@oddfoot \overfullrule 3pt
        \let\label=\draftlabel
        \let\marginnote=\draftmarginnote
   \def\@eqnnum{(\theequation)\rlap{\kern\marginparsep\tt\@eqnlabel}%
\global\let\@eqnlabel\@vacuum}  }

\def\numberbysection{\@addtoreset{equation}{section}
        \def\theequation{\thesection.\arabic{equation}}}

\def\underline#1{\relax\ifmmode\@@underline#1\else
        $\@@underline{\hbox{#1}}$\relax\fi}

\def\titlepage{\@restonecolfalse\if@twocolumn\@restonecoltrue\onecolumn
     \else \newpage \fi \thispagestyle{empty}\c@page\z@
        \def\thefootnote{\fnsymbol{footnote}} }

\def\endtitlepage{\if@restonecol\twocolumn \else  \fi
        \def\thefootnote{\arabic{footnote}}
        \setcounter{footnote}{0}}  

\hybrid

\def\beq{\begin{equation}}
\def\eq{\end{equation}}

\begin{document}

\setcounter{page}{1}


\ \vspace{0mm}

\begin{center}
 {\Large{\bf Quadratic algebras based on}}
 \\ \vspace{4mm}
 {\Large{\bf SL(NM) elliptic quantum R-matrices}}
\\
\vspace{18mm} {\large {I.A. Sechin}}\,\footnote{National Research University Higher School of Economics, Center for Advanced Studies, Skolkovo Institute of Science and Technology;
e-mail: shnbuz@gmail.com.}\qquad\quad
 {\large {A.V. Zotov}}\,\footnote{National Research University Higher School of Economics, Steklov Mathematical Institute of Russian Academy of Sciences;
  e-mail:
zotov@mi-ras.ru.}
\end{center}

\ \vskip4mm

 \begin{abstract}
    We construct quadratic quantum algebra based on the dynamical RLL--relation for
    the quantum $R$--matrix related to $SL(NM)$--bundles with nontrivial characteristic class
    over  elliptic curve.
    This $R$--matrix generalizes simultaneously the elliptic
    nondynamical Baxter--Belavin and the dynamical Felder $R$--matrices,
    and the obtained quadratic relations generalize both -- the Sklyanin
    algebra and the relations in the
    Felder--Tarasov--Varchenko elliptic quantum group,
    which are reproduced in the particular cases $M = 1$ and
    $N = 1$ respectively.
 \end{abstract}



\bigskip



\section{Sklyanin algebra}
Consider the Baxter--Belavin quantum \(R\)--matrix
\cite{Baxter, Belavin}:
\begin{equation}
  \displaystyle{
    R_{12}^{\mathrm{BB}}(\hbar, u) =
    \sum_{\alpha \in \mathbb{Z}_N^2}
        \varphi_\alpha(u, \hbar + \omega_\alpha)
        T_\alpha \otimes T_{-\alpha}\,.
 }
\end{equation}
In this definition the elliptic functions
\(\varphi_\alpha(u, x + \omega_\alpha)\) and the
\(N \times N\) basis matrices \(T_\alpha\), connected with these functions
are used. They are defined in the Appendix 1.
This \(R\)--matrix satisfies the quantum Yang--Baxter equation in
\( \mathrm{Mat}(N, \mathbb{C})^{\otimes 3} \)
\begin{equation}
    \label{YangBaxter_BB}
      \displaystyle{
    R_{12}^{\mathrm{BB}}(\hbar, z_{12})
        R_{13}^{\mathrm{BB}}(\hbar, z_{13})
            R_{23}^{\mathrm{BB}}(\hbar, z_{23}) =
        R_{23}^{\mathrm{BB}}(\hbar, z_{23})
            R_{13}^{\mathrm{BB}}(\hbar, z_{13})
                R_{12}^{\mathrm{BB}}(\hbar, z_{12})\,.
                }
\end{equation}
Here the notation \(z_{ij} = z_i - z_j \) is used, and the lower indices in
\(R\)--matrices denote the tensor component numbers, where they act nontrivially.
For instance, in (\ref{YangBaxter_BB})
\begin{equation}
  \displaystyle{
    R_{13}^{\mathrm{BB}}(\hbar, z_{13}) =
        \sum_{\alpha} \varphi_\alpha(z_{13}, \hbar + \omega_\alpha)
        T_\alpha \otimes 1_N \otimes T_{-\alpha}\,.
        }
\end{equation}
An operator \(L(z)\) is called \(L\)--operator for the
Baxter--Belavin \(R\)--matrix if it satisfies
the RLL--relation
\begin{equation}
    \label{RLL_BB}
      \displaystyle{
    R_{12}^{\mathrm{BB}}(\hbar, z_1 - z_2) L_1(z_1) L_2(z_2) =
        L_2(z_2) L_1(z_1) R_{12}^{\mathrm{BB}}(\hbar, z_1 - z_2)\,.
        }
\end{equation}
In paper \cite{Sklyanin} E. Sklyanin suggested a class of
\(L\)--operators for the case \( N = 2 \).
Later his result was extended to cases of
an arbitrary \(N\) taking also into account other possible  parameters of underlying bundles over
elliptic curve \cite{OF,Quano,CLOZ}.
The constructed \(L\)--operators are connected with the
quadratic algebra called the Sklyanin algebra.

Consider an \(L\)--operator of the form
\begin{equation}\label{u0}
  \displaystyle{
    L(z) =
    \sum_\alpha
        \varphi_\alpha(z, \hbar + \omega_\alpha)
        S_\alpha T_\alpha\,.
        }
\end{equation}
RLL--relation (\ref{RLL_BB}) for this \(L\)--operator is equivalent to
the following quadratic relations for operators \(S_\alpha\),
labelled by pairs \((\alpha, \beta)\), which do not depend on the
spectral parameters \(z_1, z_2\):
\begin{equation}\label{u1}\begin{array}{c}
  \displaystyle{
    \sum_\gamma
        \varkappa_{\gamma, \alpha-\beta}
        \Big( E_1(\omega_\gamma\! +\! \hbar) -
        E_1(\omega_{\alpha - \beta - \gamma}\! +\! \hbar) +
        E_1(\omega_{\alpha - \gamma}\! +\! \hbar) -
        E_1(\omega_{\beta + \gamma}\! +\! \hbar) \Big)
            S_{\alpha - \gamma} S_{\beta + \gamma} = 0
  }
  \end{array}\end{equation}
  for $\beta \ne 0$ and
  \begin{equation}\label{u2}\begin{array}{c}
    \displaystyle{
    \beta = 0: \quad
    \sum_\gamma
        \varkappa_{\gamma \alpha}
        \Big( E_2(\omega_\gamma + \hbar) -
        E_2(\omega_{\alpha - \gamma} + \hbar) \Big)
            S_{\alpha - \gamma} S_\gamma = 0\,,
            }
\end{array}\end{equation}
where \(E_1(z)\) and \(E_2(z)\) are the functions
defined in the Appendix 1. A set of numbers
\begin{equation}
\displaystyle{
    \varkappa_{\alpha \beta} = \exp \left(
        \frac{\pi i}{N} (\beta_1 \alpha_2 - \beta_2 \alpha_1)
    \right)
    }
\end{equation}
defines the structure constants of  relations (\ref{u1})--(\ref{u2}),
called the Sklyanin algebra relations. For example, the operators
\(S_\alpha = T_{-\alpha}\) satisfy these relations.
In this case the RLL--relation turns into the Yang--Baxter equation
(\ref{YangBaxter_BB}).

One can slightly modify the definition (\ref{u0}) and
relations (\ref{u1})--(\ref{u2}). The \(L\)--operator can be divided by a
function depending on \(z\) only, because this function is canceled in
both parts of RLL--relation. Write down \(\varphi_\alpha\)--function
explicitly
\begin{equation}
\displaystyle{
    \varphi_\alpha(z, \hbar + \omega_\alpha) =
        \phi(z, \hbar + \omega_\alpha)
            e^{\frac{2 \pi i}{N} \alpha_2 z} =
        \frac{\theta'(0) \theta(z + \hbar + \omega_\alpha)}
            {\theta(z) \theta(\hbar + \omega_\alpha)}
            e^{\frac{2 \pi i}{N} \alpha_2 z}\,
            }
\end{equation}
Dividing the \(L\)--operator (\ref{u0}) by \( \theta'(0) / \theta(z) \),
one obtains
\begin{equation}
\displaystyle{
    L^\hbar(z) =
    \sum_\alpha
        \frac{\theta(z + \hbar + \omega_\alpha)}
            {\theta(\hbar + \omega_\alpha)}
            e^{\frac{2 \pi i}{N} \alpha_2 z}
        S_\alpha T_\alpha\,.
        }
\end{equation}
The factor \(\theta(\hbar + \omega_\alpha)\) does not depend on the
spectral parameter, so that one can remove it by redefining \(S_\alpha\). In this
case the \(L\)--operator takes the form
\begin{equation}
\displaystyle{
    L^\hbar(z) =
    \sum_\alpha
        \theta(z + \hbar + \omega_\alpha)
            e^{\frac{2 \pi i}{N} \alpha_2 z}
        \widetilde{S}_\alpha T_\alpha, \quad
    \widetilde{S}_\alpha =
        \frac{S_\alpha}{\theta(\hbar + \omega_\alpha)}\,.
        }
\end{equation}
The Sklyanin algebra relations are also modified as follows:
$$\begin{array}{c}
\displaystyle{
    \beta \ne 0: \ \
    \sum_\gamma
        \varkappa_{\gamma, \alpha-\beta}
        \Big( E_1(\omega_\gamma + \hbar) -
        E_1(\omega_{\alpha - \beta - \gamma} + \hbar) +
        E_1(\omega_{\alpha - \gamma} + \hbar) -
        E_1(\omega_{\beta + \gamma} + \hbar) \Big)
         }
                \\
    \displaystyle{
         \times
        \theta(\hbar + \omega_{\alpha - \gamma})
            \theta(\hbar + \omega_{\beta + \gamma})
                \widetilde{S}_{\alpha - \gamma}
                \widetilde{S}_{\beta + \gamma} = 0\,,
                }
                 \end{array}
  $$
 \begin{equation}\label{u7}\begin{array}{c}
    \displaystyle{
    \beta = 0: \quad
    \sum_\gamma
        \varkappa_{\gamma \alpha}
        \Big( E_2(\omega_\gamma + \hbar) -
        E_2(\omega_{\alpha - \gamma} + \hbar) \Big)
        \theta(\hbar + \omega_{\alpha + \gamma})
            \theta(\hbar + \omega_\gamma)
                \widetilde{S}_{\alpha - \gamma}
                \widetilde{S}_\gamma = 0\,.}
 \end{array}\end{equation}
Moreover, one can change the parameter \(\hbar\) in the \(L\)--operators to
another parameter $\eta$ by shifting \(z\) since the \(R\)--matrix depends
 on the difference \(z_1 - z_2\) only. Then one can define the operator
\begin{equation}
\displaystyle{
    L^\eta(z) = L^\hbar(z + \eta - \hbar) =
    \sum_\alpha
        \theta(z + \eta + \omega_\alpha)
            e^{\frac{2 \pi i}{N} \alpha_2 z}
        S^\eta_\alpha T_\alpha\,, \quad
    S^\eta_\alpha = \widetilde{S}_\alpha
        e^{\frac{2 \pi i}{N} \alpha_2 (\eta - \hbar)}\,.}
\end{equation}
Relations for \(S^\eta_\alpha\) are similar to the relations for
\(\widetilde{S}_\alpha\) up to these exponential factors.

\section{Elliptic quantum group}
Consider the Felder dynamical quantum \(R\)--matrix \cite{Felder}:
\begin{equation}\label{u3}
\displaystyle{
    R_{12}^{\mathrm{F}}(\hbar, u \mid q) =
    \sum_{i = 1}^M
        \phi(u, \hbar) E_{ii} \otimes E_{ii} +
    \sum_{\substack{i, j = 1 \\ i \ne j}}^M
        \phi(u, q_{ij}) E_{ij} \otimes E_{ji} +
    \sum_{\substack{i, j = 1\\ i \ne j}}^M
        \phi(\hbar, -q_{ij}) E_{ii} \otimes E_{jj}\,,
        }
\end{equation}
where \(q_{ij} = q_i - q_j\), \(E_{ij}\) -- \(M \times M\) matrices with
matrix elements \( (E_{ij})_{kl} = \delta_{ik} \delta_{jl} \),
and \( \phi \) is the elliptic Kronecker function defined in the Appendix 1.
The title dynamical means that the \(R\)--matrix depends on the
dynamical parameters $q_i$.

The \(R\)--matrix (\ref{u3}) satisfies the quantum dynamical Yang--Baxter
equation
$$\begin{array}{c}
\displaystyle{
    R_{12}^{\mathrm{F}}(\hbar, z_{12} \mid q)
        R_{13}^{\mathrm{F}}(\hbar, z_{13} \mid q - \hbar^{(2)})
            R_{23}^{\mathrm{F}}(\hbar,z_{23} \mid q) =
            }
            \\ \ \\
            \displaystyle{
             =
        R_{23}^{\mathrm{F}}(\hbar, z_{23} \mid q - \hbar^{(1)})
            R_{13}^{\mathrm{F}}(\hbar, z_{13} \mid q)
                R_{12}^{\mathrm{F}}(\hbar, z_{12} \mid q - \hbar^{(3)})\,.
                }
\end{array}$$
In this equation the shifts along the Cartan subalgebra \(\{E_{ii}\}\) in \(\mathfrak{gl}(M)\) are used:
\begin{equation}
\displaystyle{
    R_{12}^{\mathrm{F}}(\hbar, z_{12} \mid q - \hbar^{(3)}) =
        e^{-\hbar \hat{\partial}_3}
            R_{12}^{\mathrm{F}}(\hbar, z_{12} \mid q)
                e^{\hbar \hat{\partial}_3}\,, \quad
    \hat{\partial}_3 = \sum_k (E_{kk})_3 {\partial_{q_k}}\,.
    }
\end{equation}
Besides the quantum dynamical Yang--Baxter equation,
the \(R\)--matrix satisfies also the zero--weight conditions:
$$
\displaystyle{
    [(E_{ii})_1 + (E_{ii})_2,
        R_{12}^{\rm{F}}(\hbar, z_{12} \mid q)] = 0\,,
        }
$$
$$
\displaystyle{
    [\hat{\partial}_1 + \hat{\partial}_2,
        R_{12}^{\mathrm{F}}(\hbar, z_{12} \mid q)] = 0\,.
        }
$$
Let \(h_i, i = 1, 2, \ldots M\) be commuting elements.
An operator \(L(z \mid q)\) is called the dynamical \(L\)--operator with the
Cartan elements \(h_i\) for the Felder \(R\)--matrix if it
satisfies the dynamical RLL--relation
\begin{equation}\label{u5}\begin{array}{c}
\displaystyle{
    R_{12}^{\mathrm{F}}(\hbar, z_{12} \mid q)
        L_1(z_1 \mid q - \hbar^{(2)})
            L_2(z_2 \mid q) =
        L_2(z_2 \mid q \!-\! \hbar^{(1)})
            L_1(z_1 \mid q)
                R_{12}^{\mathrm{F}}(\hbar, z_{12} \mid q \!-\! \hbar \cdot h)\,,
                }
                 \\ \ \\
                 \displaystyle{
    R_{12}^{\mathrm{F}}(\hbar, z_{12} \mid q - \hbar \cdot h) =
        e^{-\hbar \sum_k h_k \frac{\partial}{\partial q_k}}
            R_{12}^{\mathrm{F}}(\hbar, z_{12} \mid q)
                e^{\hbar \sum_k h_k \frac{\partial}{\partial q_k}}\,.}
 \end{array}\end{equation}
The dynamical Yang--Baxter equation implies that the
Felder's \(R\)--matrix is the dynamical \(L\)--operator
with the Cartan elements \(h_i = (E_{ii})_3\):
\begin{equation}
\displaystyle{
    L_1(z \mid q) = R_{13}^{\mathrm{F}}(\hbar, z \mid q)\,.
    }
\end{equation}
The RLL--relation (\ref{u5}) can be rewritten in the equivalent form
if one acts on both sides from the left by the operator
\(e^{\hbar \hat{\partial}_1} e^{\hbar \hat{\partial}_2}\).
Using the zero--weight property
\([\hat{\partial}_1 + \hat{\partial}_2, R_{12}^\hbar(u \mid q)] = 0\),
one obtains
\([e^{\hbar \hat{\partial}_1} e^{\hbar \hat{\partial}_2},
R_{12}^\hbar(u \mid q)] = 0\). Therefore, we have
$$\begin{array}{c}
\displaystyle{
    e^{\hbar \hat{\partial}_1} e^{\hbar \hat{\partial}_2}
        R_{12}^{\mathrm{F}}(\hbar, z_{12} \mid q)
            L_1(z_1 \mid q - \hbar^{(2)})
                L_2(z_2 \mid q) =
                }
                 \\ \ \\
                 \displaystyle{
                =
    e^{\hbar \hat{\partial}_1} e^{\hbar \hat{\partial}_2}
        L_2(z_2 \mid q - \hbar^{(1)})
            L_1(z_1 \mid q)
            R_{12}^{\mathrm{F}}(\hbar, z_{12} \mid q - \hbar \cdot h),
    }
 \end{array}$$
$$\begin{array}{c}
\displaystyle{
    e^{\hbar \hat{\partial}_1} e^{\hbar \hat{\partial}_2}
        R_{12}^{\mathrm{F}}(\hbar, z_{12} \mid q)
            e^{-\hbar \hat{\partial}_2}
            L_1(z_1 \mid q)
            e^{\hbar \hat{\partial}_2}
                L_2(z_2 \mid q) =
                }
                 \\ \ \\
                 \displaystyle{
                =
    e^{\hbar \hat{\partial}_2}
        L_2(z_2 \mid q)
        e^{\hbar \hat{\partial}_1}
            L_1(z_1 \mid q)
            R_{12}^{\mathrm{F}}(\hbar, z_{12} \mid q - \hbar \cdot h)\,,
            }
  \end{array}$$
  and
  $$\begin{array}{c}
  \displaystyle{
    R_{12}^{\mathrm{F}}(\hbar, z_{12} \mid q)
        e^{\hbar \hat{\partial}_1} L_1(z_1 \mid q)
            e^{\hbar \hat{\partial}_2} L_2(z_2 \mid q) =
    e^{\hbar \hat{\partial}_2} L_2(z_2 \mid q)
        e^{\hbar \hat{\partial}_1} L_1(z_1 \mid q)
        R_{12}^{\mathrm{F}}(\hbar, z_{12} \mid q - \hbar \cdot h)\,.
        }
\end{array}$$
Define the operators
\(\widetilde{L}(u \mid q) = e^{\hbar \hat{\partial}} L(u \mid q) \).
Then (\ref{u5}) can be rewritten in the form
\begin{equation}\label{u6}\begin{array}{c}
  \displaystyle{
    R_{12}^{\mathrm{F}}(\hbar, z_{12} \mid q)
        \widetilde{L}_1(z_1 \mid q)
            \widetilde{L}_2(z_2 \mid q) =
        \widetilde{L}_2(z_2 \mid q)
            \widetilde{L}_1(z_1 \mid q)
                R_{12}^{\mathrm{F}}(\hbar, z_{12} \mid q - \hbar \cdot h)\,.
                }
\end{array}\end{equation}
In paper \cite{TV} V. Tarasov and A. Varchenko constructed
the dynamical \(L\)--operators and related quadratic algebra,
which is also known as the small elliptic quantum group. Consider \(q_k\) and \(q_k - \hbar h_k\) in \(R\)--matrices in (\ref{modified_dynamic_RLL}) as independent coordinates and denote these two
new sets of variables as
\(q_k^{\{2\}} = q_k, \ q_k^{\{1\}} = q_k - \hbar h_k\).
Then the RLL--relation is written in the form
\begin{equation}
    \label{modified_dynamic_RLL}
    \begin{array}{c}
      \displaystyle{
    R_{12}^{\mathrm{F}}(\hbar, z_{12} \mid q^{\{2\}})
        \widetilde{L}_1(z_1 \mid q^{\{1\}}, q^{\{2\}})
            \widetilde{L}_2(z_2 \mid q^{\{1\}}, q^{\{2\}}) =
            }
                 \\ \ \\
                 \displaystyle{
                =
        \widetilde{L}_2(z_2 \mid q^{\{1\}}, q^{\{2\}})
            \widetilde{L}_1(z_1 \mid q^{\{1\}}, q^{\{2\}})
                R_{12}^{\mathrm{F}}(\hbar, z_{12} \mid q^{\{1\}})\,.
                }
                \end{array}
\end{equation}
Consider the following  ansatz for \(L\)--operator:
\begin{equation}
\displaystyle{
    \widetilde{L}(z \mid q) =
        \sum_{i, j} \theta(z + q_i^{\{2\}} - q_j^{\{1\}})
            t_{ji} E_{ij}\,,
            }
\end{equation}
where \(t_{ij}\) are operators, which do not commute with coordinates \(q_k^{\{I\}}\),
but shift them by \(\hbar\) according to the rule
$$\begin{array}{c}
\displaystyle{
    t_{ij}
        f(q_1^{\{1\}}, \ldots, q_i^{\{1\}}, \ldots, q_M^{\{1\}},
          q_1^{\{2\}}, \ldots, q_j^{\{2\}}, \ldots, q_M^{\{2\}})=
          }
                 \\ \ \\
                 \displaystyle{
                =
    f(q_1^{\{1\}}, \ldots, q_i^{\{1\}} + \hbar, \ldots, q_M^{\{1\}},
      q_1^{\{2\}}, \ldots, q_j^{\{2\}} + \hbar, \ldots, q_M^{\{2\}})
        t_{ij}\,,
        }
\end{array}$$
where \(f\) is an arbitrary function of variables \(q_k^{\{I\}}\).
The dynamical RLL--relation (\ref{modified_dynamic_RLL}) for this
\(L\)--operator is equivalent to the following quadratic relations
for the operators \(t_{ij}\):
$$\begin{array}{c}
\displaystyle{
    t_{ij} t_{ik} = t_{ik} t_{ij}\,,
     }
     \\ \ \\
     \displaystyle{
    t_{ik} t_{jk} =
        \frac{\theta(q_{ij}^{\{1\}} - \hbar)}
            {\theta(q_{ij}^{\{1\}} + \hbar)} t_{jk} t_{ik}\,,
        \quad i \ne j\,,
         }
        \\ \ \\
        \displaystyle{
    \frac{\theta(q_{jl}^{\{2\}} - \hbar)}
        {\theta(q_{jl}^{\{2\}})} t_{ij} t_{kl} -
    \frac{\theta(q_{ik}^{\{1\}} - \hbar)}
        {\theta(q_{ik}^{\{1\}})} t_{kl} t_{ij} =
    - \frac{\theta(\hbar) \theta(q_{ik}^{\{1\}} + q_{jl}^{\{2\}})}
        {\theta(q_{ik}^{\{1\}}) \theta(q_{jl}^{\{2\}})} t_{il} t_{kj},
        \quad i \ne k, \ j \ne l\,.
        }
\end{array}$$
These quadratic relations define the (small) elliptic
Felder--Tarasov--Varchenko quantum group.

\section{A quadratic algebra for the \(SL(NM)\) \(R\)--matrix}

Consider the quantum \(R\)--matrix related to
\(SL(NM)\)--bundle with nontrivial characteristic class over the
elliptic curve. This \(R\)--matrix was constructed in  \cite{LOSZ}.
It generalizes simultaneously the nondynamical Baxter--Belavin quantum
\(R\)--matrix and the dynamical Felder quantum \(R\)--matrix,
and can be represented in the form
\begin{equation}
\begin{array}{c}
\displaystyle{
    \mathbf{R}^\hbar_{ab12}(z_{12} \mid q)
    =\sum_i
            (E_{ii})_a (E_{ii})_b
                R_{12}^{\mathrm{BB}}(\hbar, z_{12}) +
    }
        \\ \ \\
        \displaystyle{
        +\sum_{\substack{i, j \\ i \ne j}}
            (E_{ij})_a (E_{ji})_b
                R_{12}^{\mathrm{BB}}(q_{ij}, z_{12}) +
        \sum_{\substack{i, j \\ i \ne j}}
            (E_{ii})_a (E_{jj})_b\otimes 1_N\otimes 1_N
                \phi(\hbar, -q_{ij})\,.
                }
   \end{array}
\end{equation}
Here the spaces labelled by small Latin letters are \(M \times M\) matrix spaces
in the standard basis, and the spaces labelled by numbers --- \(N \times N\) matrix
spaces in the basis (\ref{u8}). This quantum \(R\)--matrix
satisfies the dynamical quantum Yang--Baxter equation with shifts
 along the Cartan subalgebra corresponding to \(M \times M\) matrices only
(i.e. of the form $h_i\otimes 1_N$):
$$\begin{array}{c}
\displaystyle{
    \mathbf{R}^\hbar_{ab12}(z_{12} \mid q)
        \mathbf{R}^\hbar_{ac13}(z_{13} \mid q - \hbar^{(b)})
            \mathbf{R}^\hbar_{bc23}(z_{23} \mid q) =
            }
        \\ \ \\
        \displaystyle{
             =
    \mathbf{R}^\hbar_{bc23}(z_{23} \mid q - \hbar^{(a)})
        \mathbf{R}^\hbar_{ac13}(z_{13} \mid q)
            \mathbf{R}^\hbar_{ab12}(z_{12} \mid q - \hbar^{(c)})\,.
            }
\end{array}$$
An operator \( \mathbf{L}_{a1}(z \mid q^{\{1\}}, q^{\{2\}}) \) is called
an \(L\)--operator for this quantum \(R\)--matrix, if it satisfies
the following RLL--relation:
$$\begin{array}{c}
 \displaystyle{
    \mathbf{R}^\hbar_{ab12}(z_{12} \mid q^{\{2\}})
        \mathbf{L}_{a1}(z_1 \mid q^{\{1\}}, q^{\{2\}})
            \mathbf{L}_{b2}(z_2 \mid q^{\{1\}}, q^{\{2\}})
             = }
        \\ \ \\
        \displaystyle{
             =
        \mathbf{L}_{b2}(z_2 \mid q^{\{1\}}, q^{\{2\}})
            \mathbf{L}_{a1}(z_1 \mid q^{\{1\}}, q^{\{2\}})
                \mathbf{R}^\hbar_{ab12}(z_{12} \mid q^{\{1\}})\,.
                }
\end{array}$$
The main result of this paper is the description of quadratic algebra
connected with this RLL--relation. Choose an \(L\)--operator in the form
$$\begin{array}{c}
\displaystyle{
    \mathbf{L}_{a1}(z_1 \mid q^{\{1\}}, q^{\{2\}}) =
        \sum_{ij} (E_{ij})_a L^{ij}_1(z_1 \mid q^{\{1\}}, q^{\{2\}})\,, }
        \\
        \displaystyle{
    L^{ij}(z \mid q) =
        \sum_\alpha
            \theta(z + q_i^{\{2\}} - q_j^{\{1\}} + \omega_\alpha) t_{ji}^\alpha T_\alpha\,.
            }
\end{array}$$
The operators \(t_{ij}^\alpha\) shift coordinates \(q_k\) by the rule
$$\begin{array}{c}
\displaystyle{
    t_{ij}^\alpha
        f(q_1^{\{1\}}, \ldots, q_i^{\{1\}}, \ldots, q_M^{\{1\}},
          q_1^{\{2\}}, \ldots, q_j^{\{2\}}, \ldots, q_M^{\{2\}})
          =
          }
        \\ \ \\
        \displaystyle{
        =
    f(q_1^{\{1\}}, \ldots, q_i^{\{1\}} + \hbar, \ldots, q_M^{\{1\}},
      q_1^{\{2\}}, \ldots, q_j^{\{2\}} + \hbar, \ldots, q_M^{\{2\}})
        t_{ij}^\alpha\,.
        }
\end{array}$$
Then the RLL--relation is equivalent to the following set of quadratic relations
for the generators \(t_{ij}^\alpha\):
\begin{enumerate}
    \item For the same pairs of indices \( i, j \) the elements
        \( \{ t_{ji}^\alpha \mid \alpha \in \mathbb{Z}_N^2 \} \)
        satisfy the Sklyanin algebra relations with parameter
        \( \eta = q_i^{\{2\}} - q_j^{\{1\}} \).
    \item For the same second index and distinct first indices
        \(i, j, k: \ j \ne k\)
        \begin{equation}
        \displaystyle{
        \sum_\gamma
            \varkappa_{\gamma \alpha} \varkappa_{\beta \gamma}
        \phi(\hbar + \omega_\gamma,
                q_{jk}^{\{1\}} + \omega_{\beta + \gamma - \alpha})
            t_{ji}^{\alpha - \gamma} t_{ki}^{\beta + \gamma} =
        \phi(\hbar, -q_{jk}^{\{1\}}) t_{ki}^\beta t_{ji}^\alpha\,.
        }
        \end{equation}
    \item For the same first index and distinct second indices
        \(i, j, k: \ j \ne k\)
        \begin{equation}
        \displaystyle{
        \sum_\gamma
            \varkappa_{\gamma \alpha} \varkappa_{\beta \gamma}
        \phi(\hbar + \omega_{\alpha - \beta - \gamma},
                -q_{jk}^{\{2\}} - \omega_\gamma)
            t_{ik}^{\alpha - \gamma} t_{ij}^{\beta + \gamma} =
        \phi(\hbar, -q_{jk}^{\{2\}}) t_{ij}^\alpha t_{ik}^\beta\,.
        }
        \end{equation}
    \item For distinct first and second pairs of indices
        \(i, j, k, l: \ i \ne k, j \ne l\)
        \begin{equation}
        \begin{array}{c}
        \displaystyle{
        \sum_\gamma
            \varkappa_{\gamma \alpha} \varkappa_{\beta \gamma}
        \phi(q_{ik}^{\{2\}} + \omega_\gamma,
                q_{jl}^{\{1\}} + \omega_{\beta + \gamma - \alpha})
            t_{jk}^{\alpha - \gamma} t_{li}^{\beta + \gamma} =
             }
        \\
        \displaystyle{
        =
        \phi(\hbar, -q_{jl}^{\{1\}}) t_{lk}^\beta t_{ji}^\alpha -
            \phi(\hbar, -q_{ik}^{\{2\}}) t_{ji}^\alpha t_{lk}^\beta\,.
            }
            \end{array}
        \end{equation}
\end{enumerate}
In the case \(M = 1\) there are only
\(\{t_{11}^\alpha \mid \alpha \in \mathbb{Z}_N^2 \}\) generators,
satisfying the Sklyanin algebra relations,
and in the case \(N = 1\) there are only elliptic quantum groups generators
\(\{t_{ij}^0 \mid i, j \in 1, 2, \ldots, M\}\). Therefore, the constructed
quadratic algebra generalizes these two quantum algebras simultaneously.

The proof of this equivalence is straightforward, it can be verified via
elliptic functions identities given in the Appendix 1. An example of this
check for particular tensor component of the RLL--relation is presented in
the Appendix 2.

\section{Conclusion}
The quadratic algebra generalizing the elliptic quantum group and the
Sklyanin algebra is constructed. On the one hand, it is a classification type result,
which complements and generalizes the known structures of quadratic algebras
related to bundles over elliptic curve. On the other hand, the obtained results can
be applied to  description of  concrete mechanical systems. It was shown in
\cite{SeZ}, that the considered \(SL(NM)\) quantum \(R\)--matrix is connected with
quantum long--range spin chains and \(R\)--matrix--valued Lax pairs. Moreover,
this particular \(R\)--matrix in the nonrelativistic classical limit describes the
system of interacting tops. The relativistic analogue of this system was also
obtained recently using a natural ansatz for the Lax pair \cite{Z19}. So,
the result of this paper can be also considered as the description of the
operator algebra underlying the model of quantum relativistic interacting tops.

\section{Appendix}
\setcounter{equation}{0}
 \def\theequation{A.\arabic{equation}}

\subsection{Elliptic functions and their properties}
The definitions of \(R\)--matrices in this paper use the Kronecker elliptic functions
\begin{equation}\begin{array}{c}
\displaystyle{
    \varphi_\alpha(u, x + \omega_\alpha) =
        \phi(u, x + \omega_\alpha)
            e^{\frac{2 \pi i}{N} \alpha_2 u}\,, \quad
        \omega_\alpha = \frac{\alpha_1 + \alpha_2 \tau}{N},
        }
        \\ \ \\
        \displaystyle{
    \phi(u, x) = \frac{\theta'(0) \theta(u + x)}{\theta(u) \theta(x)}\,,
    }
\end{array}\end{equation}
which are expressed through the odd theta--function
\begin{equation}
\displaystyle{
    \theta(u) = -\sum_{k \in \mathbb{Z}}
        \exp \left(
            \pi i \tau \left( k + \tfrac{1}{2} \right)^2 +
            2 \pi i \left(k + \tfrac{1}{2} \right)
                \left(u + \tfrac{1}{2} \right)
        \right)\,.
        }
\end{equation}
Here \( \tau\) --- moduli of the elliptic curve, a complex parameter with \( \mathrm{Im} \tau > 0 \).

The main tool for derivation of the quadratic relations is
the addition formula (also known as the genus one Fay identity)
for the Kronecker functions
\begin{equation}
\displaystyle{
    \phi(z, x) \phi(w, y) =
        \phi(z - w, x) \phi(w, x + y) +
        \phi(w - z, y) \phi(z, x + y)
        }
\end{equation}
and its degenerations corresponding to coinciding values of variables
\begin{equation}\begin{array}{c}
\displaystyle{
    \phi(z, x) \phi(z, y) =
        \phi(z, x + y)
        (E_1(z) + E_1(x) + E_1(y) - E_1(x + y + z))\,,
        }
        \\ \ \\
        \displaystyle{
    \phi(z, x) \phi(z, -x) =
        E_2(z) - E_2(x)\,,
        }
\end{array}\end{equation}
where the Eisenstein functions \(E_1(z)\) and \(E_2(z)\) are used:
\begin{equation}\begin{array}{c}
 \displaystyle{
    E_1(z) = \frac{\theta'(z)}{\theta(z)}\,,\qquad
    E_2(z) = -E_1'(z)\,.
    }
\end{array}\end{equation}
In the definition of the Baxter--Belavin quantum \(R\)--matrix the
basis matrices \(T_\alpha\) are used. They are defined as
\begin{equation}\label{u8}\begin{array}{c}
\displaystyle{
    T_\alpha = T_{(\alpha_1, \alpha_2)} =
        \exp \left( \frac{\pi i \alpha_1 \alpha_2}{N} \right)
            Q^{\alpha_1} \Lambda^{\alpha_2}\,,
             }
        \\ \ \\
        \displaystyle{
    Q_{jk} = \delta_{jk}
        \exp \left( \frac{2 \pi i k}{N} \right)\,,
         }
\qquad
    \Lambda_{jk} = \begin{cases}
        1, \quad \text{if } j + 1 = k \mod N, \\
        0, \quad \text{else.}
    \end{cases}
\end{array}\end{equation}

\subsection{An example of calculation verifying RLL--relation}

Consider, for example, the \( (E_{ij})_a (E_{ik})_b \)--component of the
RLL--relation, for \(j \ne k\):
\begin{equation}
\begin{array}{c}
\displaystyle{
    R_{12}^{\mathrm{BB}}(\hbar, z_{12})
    L_1^{ij}(z_1) L_2^{ik}(z_2) =
    L_2^{ik}(z_2) L_1^{ij}(z_1) \phi(\hbar, -q_{jk}^{\{1\}}) +
    L_2^{ij}(z_2) L_1^{ik}(z_1)
    R_{12}^{\mathrm{BB}}(q^{\{1\}}_{kj}, z_{12})\,.
    }
\end{array}
\end{equation}
This relation is given in \(N \times N\)--matrices. Expanding it in the basis
\(T_\alpha\), one gets the following scalar relations in components \( (T_\alpha)_1 (T_\beta)_2 \)
(after cancelling all exponential factors):
$$\begin{array}{c}
\displaystyle{
    \theta(z_2 + q_i^{\{2\}} - q_k^{\{1\}} + \omega_\beta) t_{ki}^\beta \cdot
        \theta(z_1 + q_i^{\{2\}} - q_j^{\{1\}} + \omega_\alpha) t_{ji}^\alpha \cdot
        \phi(\hbar, -q_{jk}^{\{1\}}) =
         }
        \\ \ \\
        \displaystyle{
         =
    \sum_\gamma
        \varkappa_{\gamma \alpha} \varkappa_{\beta \gamma} \Big(
    \phi(z_{12}, \hbar + \omega_\gamma)
        \theta(z_1 + q_i^{\{2\}} - q_j^{\{1\}} + \omega_{\alpha - \gamma})
                t_{ji}^{\alpha - \gamma}
            \theta(z_2 + q_i^{\{2\}} - q_k^{\{1\}} + \omega_{\beta + \gamma})
                t_{ki}^{\beta + \gamma}
                 }
        \\
        \displaystyle{
                 -
    \theta(z_2 + q_i^{\{2\}} - q_j^{\{1\}} + \omega_{\alpha - \gamma})
            t_{ji}^{\alpha - \gamma}
        \theta(z_1 + q_i^{\{2\}} - q_k^{\{1\}} + \omega_{\beta + \gamma})
            t_{ki}^{\beta + \gamma}
        \phi(z_{12}, q_{kj}^{\{1\}} + \omega_{\alpha - \beta - \gamma})
    \Big)\,.
    }
\end{array}$$
Moving all \(t_{ab}\) to the right, one obtains
$$\begin{array}{c}
 \displaystyle{
    \theta(z_2 + q_i^{\{2\}} - q_k^{\{1\}} + \omega_\beta)
        \theta(z_1 + q_i^{\{2\}} - q_j^{\{1\}} + \hbar + \omega_\alpha)
        \phi(\hbar, -q_{jk}^{\{1\}}) t_{ki}^\beta t_{ji}^\alpha =
        }
        \\ \ \\
        \displaystyle{
         =
    \sum_\gamma
        \varkappa_{\gamma \alpha} \varkappa_{\beta \gamma} \Big(
    \phi(z_{12}, \hbar + \omega_\gamma)
        \theta(z_1 + q_i^{\{2\}} - q_j^{\{1\}} + \omega_{\alpha - \gamma})
            \theta(z_2 + q_i^{\{2\}} - q_k^{\{1\}} + \hbar + \omega_{\beta + \gamma})
     -
     }
        \\
        \displaystyle{
      -
    \theta(z_2 + q_i^{\{2\}} - q_j^{\{1\}} + \omega_{\alpha - \gamma})
        \theta(z_1 + q_i^{\{2\}} - q_k^{\{1\}} + \hbar + \omega_{\beta + \gamma})
            \phi(z_{12}, q_{kj}^{\{1\}} + \omega_{\alpha - \beta - \gamma})
    \Big) t_{ji}^{\alpha - \gamma} t_{ki}^{\beta + \gamma}\,.
    }
\end{array}$$
Divide two parts by
\(\theta(z_2 + q_i^{\{2\}} - q_k^{\{1\}} + \omega_\beta)
        \theta(z_1 + q_i^{\{2\}} - q_j^{\{1\}} + \hbar + \omega_\alpha)\)
and consider an expression in the brackets in the right hand side. One
can simplify it:
$$\begin{array}{c}
\displaystyle{
    \phi(z_{12}, \hbar + \omega_\gamma)
        \frac{\theta(z_1 + q_i^{\{2\}} - q_j^{\{1\}} + \omega_{\alpha - \gamma})
            \theta(z_2 + q_i^{\{2\}} - q_k^{\{1\}} + \hbar + \omega_{\beta + \gamma})}
            {\theta(z_1 + q_i^{\{2\}} - q_j^{\{1\}} + \hbar + \omega_\alpha)
        \theta(z_2 + q_i^{\{2\}} - q_k^{\{1\}} + \omega_\beta)}
     -
     }
        \\
        \displaystyle{
      -
    \phi(z_{12}, q_{kj}^{\{1\}} + \omega_{\alpha - \beta - \gamma})
        \frac{\theta(z_1 + q_i^{\{2\}} - q_k^{\{1\}} + \hbar + \omega_{\beta + \gamma})
            \theta(z_2 + q_i^{\{2\}} - q_j^{\{1\}} + \omega_{\alpha - \gamma})}
            {\theta(z_1 + q_i^{\{2\}} - q_j^{\{1\}} + \hbar + \omega_\alpha)
        \theta(z_2 + q_i^{\{2\}} - q_k^{\{1\}} + \omega_\beta)}
    =
        }
\end{array}$$
$$\begin{array}{c}
\displaystyle{
     =
    \phi(z_{12}, \hbar + \omega_\gamma)
        \frac{\phi(z_2 + q_i^{\{2\}} - q_k^{\{1\}} + \omega_\beta,
                \hbar + \omega_\gamma)}
            {\phi(z_1 + q_i^{\{2\}} - q_j^{\{1\}} + \omega_{\alpha - \gamma},
                \hbar + \omega_\gamma)}
    -
    }
        \\
        \displaystyle{
     -
    \phi(z_{12}, q_{kj}^{\{1\}} + \omega_{\alpha - \beta - \gamma})
        \frac{\phi(z_1 + q_i^{\{2\}} - q_j^{\{1\}} + \hbar + \omega_\alpha,
                q_{jk}^{\{1\}} + \omega_{\beta + \gamma - \alpha})}
            {\phi(z_2 + q_i^{\{2\}} - q_j^{\{1\}} + \omega_{\alpha - \gamma},
                q_{jk}^{\{1\}} + \omega_{\beta + \gamma - \alpha})}
    =
     }
\end{array}$$
$$\begin{array}{c}
\displaystyle{
     =
    \frac{\phi(z_{12}, \hbar + \omega_\gamma)
        \phi(z_2 + q_i^{\{2\}} - q_k^{\{1\}} + \omega_\beta,
                \hbar + \omega_\gamma)
        \phi(z_2 + q_i^{\{2\}} - q_j^{\{1\}} + \omega_{\alpha - \gamma},
                q_{jk}^{\{1\}} + \omega_{\beta + \gamma - \alpha})}
            {\phi(z_1 + q_i^{\{2\}} - q_j^{\{1\}} + \omega_{\alpha - \gamma},
                \hbar + \omega_\gamma)
            \phi(z_2 + q_i^{\{2\}} - q_j^{\{1\}} + \omega_{\alpha - \gamma},
                q_{jk}^{\{1\}} + \omega_{\beta + \gamma - \alpha})}
    -
    }
        \\ \ \\
        \displaystyle{
     -
    \frac{\phi(z_{12}, q_{kj}^{\{1\}} \!+\! \omega_{\alpha - \beta - \gamma})
        \phi(z_1 \!+\! q_i^{\{2\}} - q_j^{\{1\}} \!+\! \hbar \!+\! \omega_\alpha,
                q_{jk}^{\{1\}} \!+\! \omega_{\beta \!+\! \gamma - \alpha})
        \phi(z_1 \!+\! q_i^{\{2\}} - q_j^{\{1\}} \!+\! \omega_{\alpha - \gamma},
                \hbar \!+\! \omega_\gamma)}
            {\phi(z_1 \!+\! q_i^{\{2\}} - q_j^{\{1\}} \!+\! \omega_{\alpha - \gamma},
                \hbar \!+\! \omega_\gamma)
            \phi(z_2 \!+\! q_i^{\{2\}} - q_j^{\{1\}} \!+\! \omega_{\alpha - \gamma},
                q_{jk}^{\{1\}} \!+\! \omega_{\beta \!+\! \gamma - \alpha})}\,.
                }
\end{array}$$
Applying the Fay identity to \(\phi\), and using the property \(\phi(x, -x) = 0\),
one obtains
$$\begin{array}{c}
\displaystyle{
    \phi(z_2 + q_i^{\{2\}} - q_k^{\{1\}} + \omega_\beta,
                \hbar + \omega_\gamma)
        \phi(z_2 + q_i^{\{2\}} - q_j^{\{1\}} + \omega_{\alpha - \gamma},
                q_{jk}^{\{1\}} + \omega_{\beta + \gamma - \alpha}) =
                   }
        \\ \ \\
        \displaystyle{
                 =
    \phi(q_{jk}^{\{1\}} + \omega_{\beta + \gamma - \alpha},
            \hbar + \omega_\gamma)
        \phi(z_2 + q_i^{\{2\}} - q_j^{\{1\}} + \omega_{\alpha - \gamma},
            q_{jk}^{\{1\}} + \hbar + \omega_{\beta + 2 \gamma - \alpha})\,,
            }
        \\ \ \\
        \displaystyle{
    \phi(z_1 + q_i^{\{2\}} - q_j^{\{1\}} + \hbar + \omega_\alpha,
                q_{jk}^{\{1\}} + \omega_{\beta + \gamma - \alpha})
        \phi(z_1 + q_i^{\{2\}} - q_j^{\{1\}} + \omega_{\alpha - \gamma},
                \hbar + \omega_\gamma) =
                   }
        \\ \ \\
        \displaystyle{
                =
    \phi(\hbar + \omega_\gamma,
            q_{jk}^{\{1\}} + \omega_{\beta + \gamma - \alpha})
        \phi(z_1 + q_i^{\{2\}} - q_j^{\{1\}} + \omega_{\alpha - \gamma},
            q_{jk}^{\{1\}} + \hbar + \omega_{\beta + 2 \gamma - \alpha})\,.
            }
\end{array}$$
One can take factor out
\(\phi(\hbar + \omega_\gamma,
            q_{jk}^{\{1\}} + \omega_{\beta + \gamma - \alpha})\)  in the numerator, while
the rest parts in the numerator are equal to the denominator (through the
Fay identity):
$$\begin{array}{c}
\displaystyle{
    \phi(z_{12}, \hbar + \omega_\gamma)
        \phi(z_2 + q_i^{\{2\}} - q_j^{\{1\}} + \omega_{\alpha - \gamma},
            q_{jk}^{\{1\}} + \hbar + \omega_{\beta + 2 \gamma - \alpha})
    -
    }
        \\ \ \\
        \displaystyle{
     -
    \phi(z_{12}, q_{kj}^{\{1\}} + \omega_{\alpha - \beta - \gamma})
        \phi(z_1 + q_i^{\{2\}} - q_j^{\{1\}} + \omega_{\alpha - \gamma},
            q_{jk}^{\{1\}} + \hbar + \omega_{\beta + 2 \gamma - \alpha})
    =
    }
        \\ \ \\
        \displaystyle{
     =
    \phi(z_1 + q_i^{\{2\}} - q_j^{\{1\}} + \omega_{\alpha - \gamma},
            \hbar + \omega_\gamma)
        \phi(z_2 + q_i^{\{2\}} - q_j^{\{1\}} + \omega_{\alpha - \gamma},
            q_{jk}^{\{1\}} + \omega_{\beta + \gamma - \alpha})\,.
            }
\end{array}$$
Using this simplification, one obtains the required relation without
spectral parameters:
\begin{equation*}
\displaystyle{
    \sum_\gamma
        \varkappa_{\gamma \alpha} \varkappa_{\beta \gamma}
    \phi(\hbar + \omega_\gamma,
            q_{jk}^{\{1\}} + \omega_{\beta + \gamma - \alpha})
        t_{ji}^{\alpha - \gamma} t_{ki}^{\beta + \gamma} =
    \phi(\hbar, -q_{jk}^{\{1\}}) t_{ki}^\beta t_{ji}^\alpha\,.
    }
\end{equation*}
All other relations can be verified in the same way by considering
the other components of the RLL--relation.

\paragraph{Acknowledgements.}
%
%
This research is supported by a grant from the Russian Science Foundation
(Project No. 21-41-09011).

\begin{small}

\end{small}

\end{document}